\documentclass[11pt,reqno,english]{amsart}

\usepackage{amsmath,amsfonts,amsthm,amssymb,amsxtra}

\newtheorem{theorem}{Theorem}[section]

\theoremstyle{definition}

\theoremstyle{remark}

\newtheorem{remark}[theorem]{Remark}

\numberwithin{equation}{section}

\renewcommand{\epsilon}{\varepsilon}

\newcommand{\I}{\mathbb{I}}

\newcommand{\N}{\mathbb{N}}

\renewcommand{\phi}{\varphi}

\newcommand{\R}{\mathbb{R}}

\newcommand{\D}{-\Delta_\Omega^D}
\newcommand{\Om}{\Omega}
\newcommand{\La}{\Lambda}
\newcommand{\no}{n(\Om,\La)}

\DeclareMathOperator{\tr}{tr}
\DeclareMathOperator{\card}{card}
\DeclareMathOperator{\vol}{vol}

\title[Improved Berezin-Li-Yau inequalities with a remainder term]{Improved Berezin-Li-Yau inequalities with a remainder term}

\author{Timo Weidl}
\address{Universit\"at Stuttgart, FB Mathematik, Pfaffenwaldring 57, 70569 Stuttgart}
\email{weidl@mathematik.uni-stuttgart.de}

\begin{document}

\begin{abstract}
We give an improvement of sharp Berezin type bounds on the Riesz means $\sum_k(\Lambda-\lambda_k)_+^\sigma$
of the eigenvalues $\lambda_k$ of the Dirichlet Laplacian in a domain if $\sigma\geq 3/2$.
It contains a correction term of the order of the standard second term in the Weyl asymptotics. 
The result is based on an application of 
sharp Lieb-Thirring inequalities
with operator valued potential to spectral estimates of the Dirichlet Laplacian in domains.
\end{abstract}

\thanks{\copyright\, 2007 by the author. This paper may be reproduced, in its entirety, for non-commercial purposes.}

\dedicatory{Dedicated to M. Sh. Birman on the occasion of his 80th birthday}

\maketitle


\section{Introduction}

Let $\Om$ be an open domain in $\R^d$ and consider the Dirichlet Laplacian $\D$
defined in the form sense on the form domain $\stackrel{\circ}{H^1}\hspace{-1mm}(\Omega)$. If the embedding 
$\stackrel{\circ}{H^1}\hspace{-1mm}(\Omega)\hookrightarrow L^2(\Omega)$ is compact, for example,
if the domain $\Omega$ is bounded,
the spectrum of $-\Delta$ is discrete. It consists of a
non-decreasing sequence of positive eigenvalues $\lambda_k$ accumulating at infinity, 
which we repeat according to their  multiplicities. Let 
$$\no=\card\{\lambda_k
|\lambda_k<\La\}$$
be the corresponding counting function.

In 1912 Hermann Weyl proved the famous asymptotic formula \cite{MR1511670}
\begin{equation}\label{weylas1}
\no=(1+o(1))\eta(\Om,\La)\qquad\mbox{as}\qquad \La\to +\infty\,.
\end{equation}
where $\eta(\Omega,\Lambda)$ denotes the respective classical phase space volume
\begin{equation}\label{weylas2}
\eta(\Omega,\Lambda):=\int_{x\in\Omega}\int_{\xi\in {\mathbb
    R}^d:|\xi|^2<\Lambda}\frac{dx\cdot d\xi}{(2\pi)^d}
= L_{0,d}^{cl}{\vol}(\Omega)\Lambda^{d/2},
\quad     L_{0,d}^{cl}=\frac{\omega_d}{(2\pi)^d}\,.
\end{equation}
Here $\omega_d$ stands for the volume of the unit ball in $\R^d$.

Along side with the counting function we shall discuss the eigenvalue means
\begin{eqnarray}\label{S}
S_{\sigma,d}(\Omega,\Lambda)&:=&
\sum_n (\Lambda-\lambda_n)^{\sigma}_+
=\sigma\int_0^\Lambda (\Lambda-\tau)^{\sigma-1}n(\Omega,\tau)d\tau
\, ,
\quad \Lambda\geq 0,\ \sigma>0,\\
\label{s}
s_{\sigma,d}(\Omega,N)&:=&
\sum_{k=1}^N \lambda_k^{\sigma}
= 
\sigma\int_0^\infty \tau^{\sigma-1}(N-n(\Omega,\tau))_+d\tau
\, ,
\quad \sigma>0.
\end{eqnarray}
Inserting \eqref{weylas1} and \eqref{weylas2} into the integral \eqref{S}
one obtains the Weyl asymptotics
\begin{equation}\label{15}
S_{\sigma,d}(\Omega,\Lambda)=(1+o(1))S_{\sigma,d}^{cl}(\Omega,\Lambda)
\quad\mbox{as}\quad \Lambda\to+\infty\,,
\end{equation}
where
\begin{equation}\notag
S_{\sigma,d}^{cl}(\Omega,\Lambda)=
\int_{x\in\Omega}\int_{\xi\in{\mathbb R}^d}
(\Lambda-|\xi|^2)_+^\sigma\frac{dx\cdot d\xi}{(2\pi)^d}
=L_{\sigma,d}^{cl}{\vol}(\Omega)\Lambda^{\sigma+d/2}\,.
\end{equation}
Here we make use of the well-known Lieb-Thirring constants
\begin{equation}\notag
L_{\sigma,d}^{cl}:=
\frac{\Gamma(\sigma+1)}{2^d\pi^{d/2}\Gamma(1+\sigma+d/2)}=
\sigma B\left(\sigma,1+\frac{d}{2}\right)L_{0,d}^{cl}
\,.
\end{equation}
A similar computation for the expression \eqref{s} implies
\begin{equation}\label{16}
s_{\sigma,d}(\Omega,N)=(1+o(1))
s_{\sigma,d}^{cl}(\Omega,N)
\quad\mbox{as}\quad N\to+\infty\,,
\end{equation}
with the classical quantity
\begin{equation}\notag
s_{\sigma,d}^{cl}(\Omega,N)=
\sigma\int_0^\infty \tau^{\sigma-1}\left(N-\eta(\Omega,\tau)\right)_+d\tau
=
c(\sigma,d)\left({\vol}(\Omega)\right)^{-\frac{2\sigma}{d}}
N^{1+\frac{2\sigma}{d}}\,,
\end{equation}
where
\begin{equation}\notag
c(\sigma,d):=\frac{2\sigma}{d}\left(L_{0,d}^{cl}\right)^{-\frac{2\sigma}{d}}
B\left(\frac{2\sigma}{d},2\right)=\frac{d}{2\sigma+d}\left(L_{0,d}^{cl}\right)^{-\frac{2\sigma}{d}}
\,.
\end{equation}

It is an important question, whether the semi-classical quantities in the Weyl type formulae
can serve as universal bounds for the corresponding spectral quantities of the Dirichlet Laplacian as well.
In particular, one is interested in the validity of the  following inequalities
\begin{eqnarray}
\label{ber1}
\sum_k (\Lambda-\lambda_k)_+^\sigma=
S_{\sigma,d}(\Omega,\Lambda)&\leq& r(\sigma,d)S_{\sigma,d}^{cl}(\Omega,\Lambda)\,,
\quad\Lambda>0\,, \sigma\geq 0 \,,\\
\label{ber2}
\sum_{k=1}^N\lambda_k^\sigma=s_{\sigma,d}(\Omega,N)&\geq& \rho(\sigma,d)s_{\sigma,d}^{cl}(\Omega,N)\,,\quad
N\in\mathbb N\,,\sigma>0\,,
\end{eqnarray}
and in the sharp values of the constants $r(\sigma,d)$ and $\rho(\sigma,d)$ therein.
From the Weyl asymptotics \eqref{15}, \eqref{16} 
follows that at the best $r(\sigma,d)\geq 1$ and $\rho(\sigma,d)\leq 1$.

Berezin showed in \cite{MR0350504} that \eqref{ber1} holds with 
$r(\sigma,d)=1$ for all $d\in N$ and all $\sigma\geq 1$. 
If one applies  the Legendre transformation
to \eqref{ber1} for the
case $\sigma=1$, one finds that
\[
s_{1,d}(\Omega,N)\geq s_{1,d}^{cl}(\Omega,N)\,.
\]
The later result has independently been obtained by Li and Yau \cite{MR701919} by other means.
By H\"older's inequality we have
\begin{equation}\label{ly}
 s_{1,d}^{cl}(\Omega,N)\leq  s_{1,d}(\Omega,N)\leq s_{\sigma,d}^{1/\sigma}(\Omega,N) N^{1/\sigma'},
\quad \frac1\sigma+\frac1{\sigma'}=1\,.
\end{equation}
From this follows that \eqref{ber2} remains true
for $\sigma>1$ with the 
estimate
\[\rho(\sigma,d)\geq (c(1,d))^\sigma/c(\sigma,d)=\frac{2\sigma+d}{d}
\left(\frac{d}{2+d}\right)^\sigma
\] 
on the constant. 
If we pass in \eqref{ly} to the limit $\sigma\to\infty$, we obtain the bound
\[
\frac{d}{2+d}(L_{0,d}^{cl}\vol(\Omega))^{-\frac2d}N^{\frac2d}\leq \lambda_N\,,\quad N\in\N\,.
\]
This is equivalent to
\begin{equation}\notag
n(\Omega,\Lambda)\leq r(0,d)\eta(\Omega,\Lambda)\,,\quad\Lambda>0\,,
\end{equation}
with the estimate $1\leq r(0,d)\leq (1+2d^{-1})^{d/2}$ on the constant $r(0,d)$,
which corresponds to the case
$\sigma=0$ for the inequality \eqref{ber1}.
From here, for example,
one can conclude by the Lieb-Aizenman trick \cite{MR598768}
that $r(\sigma,d) \leq (1+2d^{-1})^{d/2}$ for $0\leq \sigma< 1$.

P\'olya showed in 1961 that for all tiling domains $\Om$ one has in fact $r(0,d)=1$, 
see \cite{MR0129219}. 
The famous P\'olya conjecture suggesting that $r(0,d)=1$ for general domains, remains unresolved so far.
For certain improvements in this direction see \cite{MR1365335,MR1491551}. Recently it has been
shown that P\'olya's conjecture fails in the case of a constant magnetic field \cite{FrLoWe}.

Along with his formula on the main high energy term for the counting function $\no$
Weyl conjectured that for the Dirichlet Laplacian the following two term formula holds true
\begin{equation}\label{twotermweyl}
n(\Omega,\Lambda)=
L_{0,d}^{cl}{\vol}(\Omega)\Lambda^{d/2}
-\frac{1}{4}L_{0,d-1}^{cl}|\partial\Omega|\Lambda^{(d-1)/2}
+o(\Lambda^{(d-1)/2})\quad\mbox{as}\quad\Lambda\to\infty\,.
\end{equation}
Here $|\partial\Omega|$ denotes the $d-1$-dimensional measure of the
boundary of $\Om$. This formula fails in general, but it 
holds true under certain restrictions on the geodesic flow in the domain \cite{MR1631419,MR1414899}. 

If we insert \eqref{twotermweyl} into \eqref{S} and \eqref{s} we obtain the following 
two-term formulae for
the eigenvalue means
\begin{eqnarray}\label{MS1}
S_{\sigma,d}(\Omega,\Lambda)&=&
L_{\sigma,d}^{cl}{\vol}(\Omega)\Lambda^{\sigma+d/2}
-\frac{1}{4}L_{\sigma,d-1}^{cl}|\partial\Omega|\Lambda^{\sigma+(d-1)/2}
+o(\Lambda^{\sigma+(d-1)/2}),\\
\notag
s_{\sigma,d}(\Omega,N)&=&c(\sigma,d)\left({\vol}(\Omega)\right)^{-\frac{2\sigma}{d}}
N^{1+\frac{2\sigma}{d}}\\
\label{MS2}
&+& 
\frac{L_{\sigma,d-1}^{cl}
(L_{\sigma,d}^{cl})^{-1-\frac{2\sigma-1}{d}}}{4(\frac{d-1}{2}+\sigma)}\cdot
\frac{\sigma|\partial\Omega|}{({\vol}(\Omega))^{1+\frac{2\sigma-1}{d}}}
N^{1+\frac{2\sigma-1}{d}}+o(N^{1+\frac{2\sigma-1}{d}})\,.
\end{eqnarray}
Of course, the sign of the second terms support the sharp inequalities stated above.

In view of these two-term asymptotics one can ask, if it is possible
to find universal bounds on the spectral quantities, which hold true for {\em all} values of $\Lambda$, 
contain the sharp first Weyl term {\em and} reflect the contribution
of the second order term?

Note that a straight-forward generalisation of the Berezin bound by just adding the second
asymptotic term on the right hand side
\[
S_{\sigma,d}(\Omega,\Lambda)\leq S_{\sigma,d}^{cl}(\Omega,\Lambda)-C\cdot 
|\partial\Omega|\Lambda^{\sigma+\frac{d-1}{2}}\quad\mbox{must fail.}
\]
Indeed, without further restrictions the r.h.s. in this bound can just be negative
for fixed values of $\Lambda$.
Therefore, any improvement of sharp one-term bounds must invoke a suitable replacement
for the surface volume $|\partial\Omega|$.

A first step towards this goal has been made by Melas \cite{MR1933356}. Refining 
the result by Li and Yau he found that
\begin{equation}\label{Melas}
\sum_{k=1}^N\lambda_k=s_{1,d}(\Omega,N)\geq s_{1,d}^{cl}(\Omega,N)+
M(d)\frac{{\vol}(\Omega)}{J(\Omega)}N\,.
\end{equation}
Here
\[J(\Omega)=\min_{y\in{\mathbb R}^d} \int_\Omega |x-y|^2dx\]
is the momentum of $\Omega$ and the constant $M(d)$ depends only on the dimension $d$.
This bound is remarkable, since it improves the Li-Yau inequality
which corresponds via the Legendre transformation
to the endpoint $\sigma=1$
of the scale, where $r(\sigma,d)=1$ is known to be true. On the other hand, comparing
\eqref{Melas} with \eqref{MS2} we find, that the additional term of order $N$ 
on the r.h.s. of \eqref{Melas} does not capture the
correct order $N^{1+\frac{1}{d}}$ of the second Weyl term.

In the present note we give with \eqref{mainres} in Theorem 2.1 an improvement of a Berezin type bound 
\eqref{ber1} for $\sigma\geq 3/2$.
The result is based on the application of 
sharp Lieb-Thirring inequalities
with operator valued potentials 
\cite{MR1756570,MR1775696} to spectral estimates of the Dirichlet Laplacian in domains, 
see also \cite{MR2107707,MR1883337}. In fact, we 
simply raise the explicit remainder term
for the second derivative with Dirichlet boundary condition on one-dimensional intervals
to the Dirichlet Laplacian on domains in arbitrary dimensions.
Although the proof is straightforward, we find
\eqref{mainres} noteworthy, since
\begin{itemize}
\item it contains a correction term of the order of the standard second term in the Weyl asymptotics,
\item it improves even the first term by taking only the volume of a suitable subset $\Om_\Lambda\subset\Om$ into account,
\item it is applicable even for certain domains $\Omega$ with infinite volume,
\item it generalises without change of constants to Dirichlet Laplacians with arbitrary magnetic fields.
\end{itemize}
However, our result does not touch the endpoint $\sigma=1$ of the scale where the sharp one-term
Berezin bound is known. Hence, it does not imply an improvement on the best known estimates on the
counting function. The question, whether \eqref{mainres} can be generalised in some way to $\sigma=1$
or whether the Melas bound can be improved with a correction term of the second Weyl term order, remains open.

The structure of the paper is as follows. In section 2 we state the main result and comment on
the improvements which we gain compared to the one-term bound. The proof of our result will be sketched in
section 3.

\section{Statement of the Result}

\subsection*{Notation}

First let prepare some notation which shall be of use below.

Let $\Omega$ be an open subset of $\R^d$. 
We fix a Cartesian coordinate system in ${\mathbb R}^d$ and 
for $x\in\R^d$ we shall also write 
$x=(x',t)\in {\mathbb R}^{d-1}\times\mathbb R$.
Each of the sections
$\Omega(x')=\{t\in\R|(x',t)\in\Omega\}$ consists
of at most countably many open intervals $J_k (x')\subset\R$ 
of length $l_k(x')$. 
For given $x'$ and $\Lambda$ let $\kappa(x',\Lambda)\subset\N$ be the subset of 
all indices $k$, where $l_k(x')>l_\Lambda:=\pi\Lambda^{-1/2}$. 

If $\vol(\Omega_\Lambda)$ is finite, 
the sets $\kappa(x',\Lambda)$ are finite for a.e. $x'$ and we let
$\varkappa(x',\Lambda)$ denote the number of the elements of these sets. 
We shall assume that this function is measurable in $x'$.

Put
\begin{equation}\notag
\Omega_\Lambda(x')=\bigcup_{k\in\kappa(x',\Lambda)}
J_k (x')\subset\Omega(x')\subset\R
\quad\mbox{and}\quad 
\Omega_\Lambda=\bigcup_{x'\in{\mathbb R}^{d-1}}  
\{x'\}\times\Omega_\Lambda(x')\subset
\Omega\,.
\end{equation}
That means $\Omega_\Lambda$ is the subset of $\Omega$, where the intervals  in $t$-direction
contained in
$\Omega$ are longer than $l_\Lambda$. 
The set $\Omega_\Lambda$ is monotone increasing in $\Lambda$.

We shall also make use of the quantity
\begin{equation}\notag
d_\Lambda(\Omega)=\int_{x'\in\R^{d-1}} \varkappa(x',\Lambda)dx'\,,
\end{equation}
which is an effective measure of the projection
of $\Omega_\Lambda$ on the $x'$-plane counting the number of sufficiently
long intervals.

\subsection*{A basic estimate}

Consider the function
\[
f_\mu(A)=
\frac{A}{2}B\left(1+\mu,\frac{1}{2}\right)-
\sum_{k=1}^\infty \left(1-\frac{k^2}{A^2}\right)^{\mu}_+\,,
\quad A\geq 1\,,\mu>0\,,
\]
which is continuous in $A$. Since
\[\frac{A}{2}B\left(1+\mu,\frac{1}{2}\right)=\int_0^{+\infty} \left(1-\frac{t^2}{A^2}\right)_+^\mu dt
> \sum_{k=1}^\infty \left(1-\frac{k^2}{A^2}\right)^{\mu}_+\,,\]
the function $f_\mu(A)$ takes positive values for $A\geq 1$. Moreover, it is
not difficult to see, that
\[
\lim_{A\to+\infty} f_\mu(A)=\frac12\,,\quad\mu>0\,.
\]
Thus, the minimum
\[\varepsilon_\mu=\min_{A\geq 1} f_\mu(A)>0\]
exists and
\begin{equation}\label{epsi}
\sum_{k=1}^\infty \left(1-\frac{k^2}{A^2}\right)^{\mu}_+\leq \frac{A}{2}B\left(1+\mu,\frac{1}{2}\right)-\varepsilon_\mu\,,
\quad A\geq 1\,.
\end{equation}
The cases $A=1$ and $A\to +\infty$ ensure that
\[
2\varepsilon_{\mu}
\leq \min\left\{1,B\left(1+\mu,\frac12\right)\right\}\,.
\]

\subsection*{Main Result}

\begin{theorem}\label{tm1}
There exists a positive constant $\nu=\nu(\sigma,d)$, such that
for any open domain $\Omega\subset{\mathbb R}^d$, $\sigma\geq 3/2$ 
and any $\Lambda>0$ the bound
\begin{equation}\label{mainres}
\sum_{k} (\Lambda-\lambda_k)_+^\sigma
= S_{\sigma,d}(\Omega,\Lambda)
\leq L_{\sigma,d}^{cl}{\vol}(\Omega_\Lambda)\Lambda^{\sigma+\frac{d}{2}}
-\nu(\sigma,d)
\frac{L_{\sigma,d-1}^{cl}}{4}d_\Lambda(\Omega)
\Lambda^{\sigma+\frac{d-1}{2}}
\end{equation}
holds true.
The optimal constant $\nu=\nu(\sigma,d)$ satisfies
\begin{equation}\label{4e}
4\varepsilon_{\sigma+\frac{d-1}{2}}\leq\nu(\sigma,d)
\leq 2\min\left\{1,B\left(1+\sigma+\frac{d-1}{2},\frac12\right)\right\}\,.
\end{equation}
\end{theorem}

\subsection*{Comments}

\begin{remark}
The correction term in \eqref{mainres} reflects
the correct asymptotical order $O(\Lambda^{\sigma+\frac{d-1}{2}})$ of
the standard correction term in the Weyl formula.
\end{remark}

\begin{remark}
The first term on the r.h.s. of \eqref{mainres} takes only the volume
of the part $\Omega_\Lambda\subset\Omega$ into account, that is
only the part of the original domain where it is sufficiently wide for a
Dirichlet bound state  below $\Lambda$ to settle in $t$-direction.
Therefore, we have already an improvement in the first term of the Berezin bound.
In particular, the inequality \eqref{mainres} is applicable to domains $\Om$ of infinite
volume, as long as for given $\Lambda$ the volume of $\Om_\Lambda$ is finite.
\end{remark}

\begin{remark}
Numerical evaluations which can be made rigorous by elementary means show that 
the lower bound on $\nu(\sigma,d)$ is reasonable. For example, if $d=2$ and $\sigma=3/2$,
then
\[
1.91<\nu\left(\frac32,2\right)\leq 2\,.
\]
\end{remark}

\subsection*{The upper bounds}
At this point
let us discuss the upper bounds on $\nu(\sigma,d)$ more in detail. First of all, 
consider a long, thin rectangle 
of width $w$ and height $h=\delta w$, which is parallel to the coordinate axes. 
For sufficiently large $\Lambda$ we have $\Om=\Om_\Lambda$, $d_\Lambda=w$ and 
$|\partial\Om|=2(1+\delta)d_\lambda$.
The asymptotical formula \eqref{MS1} in comparison with \eqref{mainres} yields  that 
$\nu(\sigma,d)\leq 2(1+\delta)$. Since $\delta$ can be chosen arbitrary small, one finds that
$\nu(\sigma,d)\leq 2$.

Now take again a rectangle of width $w$ and height $h=(1-\delta)l_\Lambda$ for some 
small positive
$\delta$. Then $d_\Lambda=w$, $\Om=\Om_\Lambda$ and 
$\vol(\Om_\Lambda)=(1-\delta)\pi\Lambda^{-1}d_\Lambda$. The r.h.s. of \eqref{mainres}
turns into 
\[
 d_\Lambda(\Omega)\Lambda^{\sigma+\frac{d-1}{2}}
\left((1-\delta)\pi L_{\sigma,d}^{cl} - \frac{\nu(\sigma,d)}{4}L_{\sigma,d-1}^{cl}\right)\,.
\]
Since this expression has to be non-negative for the bound \eqref{mainres} to make sense, we find that 
$$\nu(\sigma,d)\leq 4\pi  L_{\sigma,d}^{cl}(L_{\sigma,d-1}^{cl})^{-1}=
2B\left(\frac12,1+\sigma+\frac{d-1}{2}\right)\,.$$

Actually, this condition ensures that  the r.h.s. of \eqref{mainres}
is non-negative for arbitrary $\Om$. Indeed,
since the total width of $\Om_\Lambda(x')$ exceeds
$\varkappa(x',\Lambda)\cdot l_\Lambda\geq \pi\Lambda^{-1/2}\varkappa(x',\Lambda)$, the volume
of $\Omega_\Lambda$ can be estimated from below by
\[
\vol(\Omega_\Lambda)=\int_{x'\in{\mathbb R}^{d-1}} \varkappa(x',\Lambda)\cdot l_\Lambda dx'
\geq \pi\Lambda^{-1/2} d_\Lambda(\Omega).
\]
Thus, the expression on the r.h.s. of \eqref{mainres} satisfies
\[
L_{\sigma,d}^{cl}{\vol}(\Omega_\Lambda)\Lambda^{\sigma+\frac{d}{2}}
-
\frac{\nu(\sigma,d)}{4}L_{\sigma,d-1}^{cl}d_\Lambda(\Omega)
\Lambda^{\sigma+\frac{d-1}{2}}
\geq d_\Lambda(\Omega)\Lambda^{\sigma+\frac{d-1}{2}}
\left(\pi L_{\sigma,d}^{cl} - \frac{\nu(\sigma,d)}{4}L_{\sigma,d-1}^{cl}\right)\geq 0\,.
\]

\section{The proof of the main result }

The proof proceeds in two steps
We start with a variational argument, which transforms our initial problem into a
spectral estimate on the negative eigenvalues of a Schr\"odinger type operator with
operator values potential.
Our result follows then from an operator valued Lieb-Thirring bound
\cite{MR1756570,MR1760755,MR1775696}.  

\subsection*{Step 1: A variational argument}

We use the notation introduced in section 2. Moreover,
by $\nabla^\prime$ and $-\Delta^\prime$ we denote the gradient and
the Laplacian in the first $d-1$ dimensions of $\R^d=\R^{d-1}\times\R\ni(x',t)$.

Consider a function $u$ from the form core $C_0^\infty(\Om)$  of $-\Delta_\Om^D$. For the
quadratic form of $-\Delta_\Om^D-\Lambda$ we have the identity
\begin{equation*}
\|\nabla u\|^2_{L^2(\Om)}-\Lambda\|u\|^2_{L^2(\Om)} = \|\nabla^\prime u\|^2_{L^2(\Om)}
+\int_{\R^{d-1}} dx' \int_{\Omega(x')} 
\left(\left|\frac{\partial}{\partial t}u(x',t)\right|^2-\Lambda|u(x',t)|^2\right)dt\,.
\end{equation*}

Note that $\Omega(x')=\cup_k J_k(x')$.
The functions $u(x',\cdot)$ satisfy Dirichlet boundary conditions
at the endpoints of $J_k(x')$. The lowest Dirichlet eigenvalue of $-d^2/dt^2$
on $J_k(x')$ equals $\pi^2l_k^{-2}(x')$. Hence, 
\[
\int_{J_k(x')} \left|\frac{\partial}{\partial t}u(x',t)\right|^2dt\geq 
\frac{\pi^2}{l_k^2(x')}\|u(x',\cdot)\|^2_{L^2(J_k(x'))}\,.
\]
In particular, for all $k\notin \kappa(x',\Lambda)$ we have
\[
\int_{J_k(x')} 
\left(\left|\frac{\partial}{\partial t}u(x',t)\right|^2-\Lambda|u(x',t)|^2\right)dt
\geq 0\,.
\]
It follows that
\begin{eqnarray*}
\int_{\Omega(x')} 
\left(\left|\frac{\partial}{\partial t}u(x',t)\right|^2-\Lambda|u(x',t)|^2\right)dt
&\geq&  \sum_{k\in\kappa(x',\Lambda)}  \int_{J_k(x')} 
\left(\left|\frac{\partial}{\partial t}u(x',t)\right|^2-\Lambda|u(x',t)|^2\right)dt\\
&\geq&
-\sum_{k\in\kappa(x',\Lambda)} 
\left<V_ku(x',\cdot),u(x',\cdot)\right>_{L^2(J_k(x'))}
\end{eqnarray*}
where the operators $V_k=V_k(x',\Lambda)$ denote the negative parts of the Sturm-Liouville problems
$-\frac{d^2}{dt^2}-\Lambda$ with Dirichlet boundary conditions on $L^2(J_k(x'))$.
Put
\[
W(x',\Lambda)=\left(\oplus_{k\in\kappa(x',\Lambda)}V_k(x',\Lambda)\right)\oplus\mathbb{O}
\quad\mbox{on}\quad 
L^2(\Omega_\Lambda(x'))\oplus L^2(\Omega(x')\setminus\Omega_\Lambda(x'))\,.
\]
This operator is bounded on $L^2(\Omega(x'))$ and for any $u\in C_0^\infty(\Omega)$ it holds
\begin{equation*}
\|\nabla u\|^2_{L^2(\Om)}-\Lambda\|u\|^2_{L^2(\Om)} \geq \|\nabla^\prime u\|^2_{L^2(\Om)}
-\int_{\R^{d-1}}\left<W(x',\Lambda)u(x',\cdot),u(x',\cdot)\right>_{L^2(\Omega(x'))}dx'\,.
\end{equation*}
In a last step we study now the function $f=u+v$ where both functions $u\in C_0^\infty(\Omega)$ and
$v\in C_0^\infty(\hat{\Om})$ are extended by zero to $\R^d$ and 
where $\hat{\Om}=\R^d\setminus\overline{\Omega}$. Since 
$\|\nabla v\|^2_{L^2(\hat{\Om})}\geq
\|\nabla^\prime v\|^2_{L^2(\hat{\Om})}$ we have
\begin{equation}\label{var}
\|\nabla v\|^2_{L^2(\hat{\Om})}+
\|\nabla u\|^2_{L^2(\Om)}-\Lambda\|u\|^2_{L^2(\Om)} \geq \|\nabla^\prime f\|^2_{L^2(\R^d)}
-\int_{\R^{d-1}}\left<W f(x',\cdot),f(x',\cdot)\right>_{L^2(\R)}dx'\,.
\end{equation}
Here we extend (in a slight abuse of notation) 
$W=W(x',\Lambda)$ by an orthogonal sum with $\mathbb{O}|_{L_2(\R\setminus\overline{\Om(x')})}$
to a bounded operator on $L^2(\R)$. 

The inequality \eqref{var} holds true for $f\in  C_0^\infty(\R^d\setminus\partial\Omega)$,
which is a form core for $(-\Delta_{\hat{\Omega}}^D)\oplus(-\Delta_\Om^D-\Lambda)$
corresponding to the expression on the l.h.s.
The semi-bounded form on the r.h.s. is closed on the larger domain
$H^1(\R^{d-1},L^2(\R))$,
where it corresponds to the Schr\"odinger type operator
$-\Delta^\prime\otimes\I - W(x',\Lambda)$ on $L^2(\R^{d-1},L^2(\R))$. Due to the positivity of 
$-\Delta_{\hat{\Omega}}^D$ the variational principle implies that
\begin{eqnarray}
\notag
\sum_k (\Lambda-\lambda_k)^\sigma_+=\tr (-\Delta_\Om^D-\Lambda)^\sigma_-
&=&\tr \left((-\Delta_{\hat{\Omega}}^D)\oplus(-\Delta_\Om^D-\Lambda)\right)^\sigma_-\\
\label{33}
&\leq& \tr \left(-\Delta^\prime\otimes\I - W(x',\Lambda)\right)_-^\sigma\,.
\end{eqnarray}

\subsection*{Step 2: Sharp Lieb-Thirring bounds}
We are now in the position to apply a sharp Lieb-Thirring inequality for operator valued
potentials (Theorem 3.1 in \cite{MR1756570}) in the dimension $d-1$ with $\sigma\geq 3/2$. 
In our setting it reads as follows
\begin{equation}\label{32}
\tr \left(-\Delta^\prime\otimes\I - W(x',\Lambda)\right)_-^\sigma
\leq L_{\sigma, d-1}^{cl}\int_{\R^{d-1}}\tr W^{\sigma+\frac{d-1}{2}}(x',\Lambda)dx'\,.
\end{equation}

The eigenvalues of the operator $W(x',\Lambda)$ are known explicitly and the non-zero
ones equal $\mu_{j,k}=\Lambda-j^2\pi^2l_k^{-2}(x')$ for $k\in\kappa(x',\Lambda)$ and
$j=1,2,\dots , [l_k(x')l_\Lambda^{-1}]$.
If we insert this into \eqref{32} and \eqref{33} we obtain
\begin{eqnarray}\notag
\sum_k (\Lambda-\lambda_k)^\sigma_+&\leq&  L_{\sigma, d-1}^{cl}\int_{\R^{d-1}} 
\sum_{k\in\kappa(x',\Lambda)} \sum_{j=1}^{ [l_k(x')l_\Lambda^{-1}]}
(\Lambda-j^2\pi^2l_k^{-2}(x'))^{\sigma+\frac{d-1}{2}}dx'\\
\label{34}
&\leq&
\Lambda^{\sigma+\frac{d-1}{2}}L_{\sigma, d-1}^{cl}\int_{\R^{d-1}} 
\sum_{k\in\kappa(x',\Lambda)} \sum_{j=1}^{ [l_k(x')l_\Lambda^{-1}]}
\left(1-\frac{j^2}{l_k^{2}(x')l^{-2}_\Lambda}\right)^{\sigma+\frac{d-1}{2}}dx'\,.\quad
\end{eqnarray}
Since for $k\in\kappa(x',\Lambda)$ we have $l_k(x')l_\Lambda^{-1}\geq 1$. From
\eqref{epsi} it follows that
\[
\sum_{j=1}^{ [l_k(x')l_\Lambda^{-1}]} \left(1-\frac{j^2}{l_k^2(x')l_\Lambda^{-2}}\right)^{\sigma+\frac{d-1}{2}}
\leq
\frac{l_k(x')l_\Lambda^{-1}}{2}B\left(1+\sigma+\frac{d-1}{2},\frac12\right)-
\varepsilon_{\sigma+\frac{d-1}{2}}
\,,
\]
what together with \eqref{34} amounts to
\begin{eqnarray}\notag
\sum_k (\Lambda-\lambda_k)^\sigma&\leq& 
\frac{\Lambda^{\sigma+\frac{d}{2}}}{2\pi}B\left(1+\sigma+\frac{d-1}{2},\frac12\right) L_{\sigma, d-1}^{cl}\int_{\R^{d-1}}\sum_{k\in\kappa(x',\Lambda)}  l_k(x')dx'\\
\label{fin1}
&-& \varepsilon_{\sigma+\frac{d-1}{2}}
\Lambda^{\sigma+\frac{d-1}{2}}L_{\sigma, d-1}^{cl}\int_{\R^{d-1}} 
\sum_{k\in\kappa(x',\Lambda)} 1 dx'\,.
\end{eqnarray}
Note that $\sum_{k\in\kappa(x',\Lambda)}  l_k(x')$ stands for the total width of $\Omega_\Lambda$
in $t$-direction and
\begin{equation}\label{fin2}
\int_{\R^{d-1}}\sum_{k\in\kappa(x',\Lambda)}  l_k(x')dx'=\vol(\Omega_\Lambda)\,.
\end{equation}
Moreover, we have
\begin{equation}\label{fin3}
\int_{\R^{d-1}} 
\sum_{k\in\kappa(x',\Lambda)} 1 dx' = \int_{\R^{d-1}} \varkappa(x',\Lambda)dx'=d_\Lambda(\Omega)\,.
\end{equation}
A straightforward computation shows that 
\begin{equation}\label{fin4}
\frac{1}{2\pi}B\left(1+\sigma+\frac{d-1}{2},\frac12\right) L_{\sigma, d-1}^{cl}= L_{\sigma, d}^{cl}\,.
\end{equation}
Inserting \eqref{fin2}-\eqref{fin4} into \eqref{fin1}, we complete the proof.
 
\subsection*{The magnetic case}
Here we consider the magnetic Dirichlet Laplacian $(i\nabla+A(x))^2$ defined in the forms sense from 
a form core
$C_0^\infty(\Omega)$ for $A\in L^d_{loc}(\Omega,\R^d)$ for $d\geq 3$ and 
$A\in L^{1+\epsilon}_{loc}(\Omega,\R^d)$
for $d=2$ and some $\epsilon>0$. It turns out that the bound \eqref{mainres} holds with the same estimates
on $\nu(\sigma,d)$ in the magnetic case as well. 
The proof uses the idea of induction in dimension in the magnetic case,
see \cite{MR1756570} section 3.2, an argument due to Helffer. 
We sketch it for the benefit of the reader.

First assume for simplicity that the vector potential $A$ is continuous.
For $x=(x_1,\dots,x_d)$ let $x=(x'_j,x''_j)$ with $x'_j=(x_1,\dots,x_j)$ and $x''_j=(x_{j+1},\dots,x_d)$.
We put 
$$\Omega(x'_j)=\{x''_j\in\R^{d-j}|x=(x'_j,x''_j)\}\in\Omega.$$
This is an open set in $\R^{d-j}$ and let $(i\nabla''_j+A''(x''))^2$ be the magnetic Dirichlet Laplacian
on it. In a first step we can chose a gauge $A(x)=(0,A''_1(x'_1,x''_1))$, where the first component of $A(x)$ 
vanishes. Using a variational argument similar to the one 
above and a sharp Lieb-Thirring bound in the dimension one,
we find that for $\sigma\geq 3/2$
\[
\sum_k(\Lambda-\lambda_k)^\sigma_+=
\tr ((i\nabla+A(X))^2-\Lambda)_-^\sigma \leq L_{\sigma,1}^{cl} 
\int_\R\tr ((i\nabla''_1+A''_1(x_1',x_1''))^2-\Lambda)_-^{\sigma+\frac12}dx_1\,.
\]
We apply now the same argument to the operator $(i\nabla''_1+A''_1(x_1',x_1''))^2$ on $\Omega(x'_1)$
an so on. Decreasing dimensions one by one we get in the $(d-1)$st step
\[\sum_k(\Lambda-\lambda_k)^\sigma_+\leq 
 L_{\sigma,1}^{cl} L_{\sigma+\frac{1}{2},1}^{cl}\cdots  L_{\sigma+\frac{d-2}{2},}^{cl}
\int_{\R^{d-1}}\tr \left(\left(i\frac{d}{dx_d}+\tilde{A}''_d(x'_{d-1},x_d)\right)^2-\Lambda\right)_-^{\sigma+\frac{d-1}{2}}dx_1\,.
\]
With the notation $\tilde{A}''_d$ we take into account, that the explicit expression for the 
vector potential changes at each step due to the necessary gauge transformation.

We avail at a family of one-dimensional Sturm-Liouville problems, where the magnetic
field simply can be gauged away. That means we can make use of the same arguments for the
operator in $t=x''_{d-1}$-direction as above and avail at 
\begin{equation}\notag
\sum_k (\Lambda-\lambda_k)^\sigma_+
\leq
\Lambda^{\sigma+\frac{d-1}{2}}L_{\sigma, 1}^{cl}\cdots L_{\sigma+\frac{d-2}{2}, 1}^{cl}\int_{\R^{d-1}}
\sum_{k\in\kappa(x',\Lambda)} \sum_{j=1}^{ [l_k(x')l_\Lambda^{-1}]}
\left(1-\frac{j^2}{l_k^{2}(x')l^{-2}_\Lambda}\right)^{\sigma+\frac{d-1}{2}}dx'\,.\quad
\end{equation}
Since $L_{\sigma, 1}^{cl}\cdots L_{\sigma+\frac{d-2}{2}, 1}^{cl}= L_{\sigma, d-1}^{cl}$, we complete the
proof in the same way as above and find
\begin{equation}\label{mainresmg}
\tr ((i\nabla+A(x))^2-\Lambda)_-^\sigma
\leq L_{\sigma,d}^{cl}{\vol}(\Omega_\Lambda)\Lambda^{\sigma+\frac{d}{2}}
-\nu(\sigma,d)
\frac{L_{\sigma,d-1}^{cl}}{4}d_\Lambda(\Omega)
\Lambda^{\sigma+\frac{d-1}{2}}
\end{equation}

Since the estimates\footnote{This does not mean that the constants $\nu(\sigma,d)$ itself
are necessarily independent of $A$.}
 on the constants $\nu(\sigma,d)$ do not depend on the magnetic vector potential, we
can close the result to non-smooth $A$ by a standard argument. 

For this let just note that in the case of a counter example,
even if $\Lambda$ is an accumulation point of the spectrum, 
already a finite partial eigenvalue sum will fail the bound \eqref{mainresmg}. Thus,
it suffices to study the quadratic form $\int_\Om|(i\nabla+A)^2u|^2dx$ on the
finite-dimensional subspace of $u$ spanned by the corresponding eigenfunctions. 
In fact, one can reduce oneself to a finite-dimensional subspace of the form core.
Since $A$ can be approximated by continuous $A$ in $L^2_{loc}$, we have a uniform convergence
$\int_\Om|(i\nabla+A_\delta)^2u|^2dx\to\int_\Om|(i\nabla+A)^2u|^2dx$ on normalised $u$ from
this subspace. Variational arguments give a contradiction to the counter example.

This line of arguments will always works, if we have an approximation of the vector potential
by continuous ones in terms of
the convergence of the forms on any finite dimensional subspace of a form core in 
$C_0^\infty(\Omega)$. In particular, it applies to the minimal extensions of theAharonov-Bohm
type operators for $d=2$.

In conclusion we remark, that the standard one-term Berezin bound
\begin{equation}\label{mgber}
\tr ((i\nabla+A(x))^2-\Lambda)_-^\sigma
\leq L_{\sigma,d}^{cl}{\vol}(\Omega)\Lambda^{\sigma+\frac{d}{2}}
\end{equation}
holds true
with
the semi-classical constant for constant magnetic 
fields if $\sigma\geq 1$ \cite{MR1779898}. It fails even in this special case  
for $\sigma<1$ \cite{FrLoWe}.
The question, whether \eqref{mgber} holds true (or can even be improved by a 
suitable correction term) 
for general magnetic fields for
$\sigma\in[1,3/2)$ 
remains open.

\section{Acknowledgements}
The author acknowledges support by the DAAD D/04/26013 and by the DFG
We 1964-2/1. He is grateful to R. Benguria, P. Exner and A. Laptev for helpful comments.
The material has been discussed and improved during the workshop  
`Low eigenvalues of Laplace and Schr\"odinger operators'
held at AIM in May 2006.

\end{document}